\documentclass[12pt]{amsart}

\newcommand\real{{\mathrm I}\!{\mathrm R} }

\newcommand\rat{{\mathrm Q}\kern-.65em {}^{{}_/ }}

\newtheorem{corollary}{Corollary}
\newtheorem{remark}{Remark}

\newtheorem{theorem}{Theorem}
\newtheorem{lemma}{Lemma}
\newtheorem{proposition}{Proposition}

\begin{document}
\title{On the number of periodic geodesics in rank 1 surfaces}
\author{Abdelhamid Amroun}
\date{January 2016}
\maketitle
\address{\begin{center}
Laboratoire de Math\'ematiques d'Orsay, Univ. Paris-Sud, CNRS, 
Universit\'e Paris-Saclay, 91405 Orsay, France.
\end{center}}
\begin{abstract} 
We consider the geodesic flow of a compact connected rank 1 surface. 
We prove a formula for the topological pressure as the exponential growth rate of rank 1 
periodic geodesics generalizing a previous result of K. Gelfert and B. Schapira \cite{GS}.
\end{abstract}
\section{Introduction and main result}
In \cite{GS} K. Gelfert and B. Schapira, have proved a formula
for a compact connected surfaces which express the topological pressure of a class of potentials
as the exponential growth rate of rank 1 geodesics with Lyapunov exponent bounded bellow by a positive constant depending on the potential. The potential they consider are H\"older continuous, hyperbolic (in the sens of \cite{GS} Theorem 1.2 and Lemma 3) and constant on the subset of singular unit tangent vectors to the surface.
In this paper, we want to generalize this formula to all 
H\"older continuous potentials which are constant on the singular part, i.e which are not necessarily hyperbolic. This leads to a formula which includes periodic orbits close to the singular part
(Lyapunov exponents may be close to zero ).
Before we explain the details and state precisely the main result, we introduce the notations and the settings. In what follows, by curvature we mean sectional curvature.

Let $M$  be a compact connected rank 1 manifold
and $G=(g^{t})_{t\in \real}$ its geodesic flow. 
According to \cite{BuSp}, the rank 1 condition is generic in the set of
nonpositively curved compact manifolds, i.e ``most'' of the compact manifolds
of nonpositive curvature are rank 1.
The rank 1 condition implies the existence of a rank 1 geodesic, i.e 
which does not admit a parallel perpendicular Jacobi field.
The space of all parallel Jacobi fields along these geodesics is then one dimensional. 
We say that a rank 1 geodesic is hyperbolic or regular.
A tangent vector $v \in TM$ has rank 1 if it defines a rank 1 geodesic. 
The set $reg$ of all the rank 1 vectors (regular vectors) is open and dense in $T^{1}M$ \cite{Ballmann}
and invariant with respect to the action of $G$ on $T^{1}M$.
The singular part $sing$ of $T^{1}M$ consists of all non-regular 
vectors in $T^{1}M$; the curvature vanishes along any singular geodesic (a geodesic defined by a singular vector).
The existence of singular geodesics implies that the geodesic flow is not Anosov.
For a complete exposition of this subject see \cite{eber, eberl, eberll}. 

Let $G(M)$ be the set of all primitive closed (periodic) geodesics, 
which are distinct and represent different free homotopy classes. 
Note that two distinct regular closed geodesics
can not be contained in the same free homotopy class, otherwise, they 
will bound a flat strip which contradicts the regularity property.

Given $\Gamma \subset G(M)$ we set
\begin{equation*}
P_{\Gamma}(t):=\{\gamma \in \Gamma : l(\gamma)\leq t\}.
\end{equation*}
Here $l(\gamma)$ is the least period of the closed geodesic $\gamma$. 
The number of periodic geodesics in a free homotopy class may be infinite but
the corresponding periods are equal.

Set $P(t):=P_{G(M)}(t)$ and let $P_{reg}(t)\subset P(t)$ be the subset
of rank 1 geodesics. 
It follows from the result of G. Knieper \cite{knieper2} that,
\begin{equation}
\lim_{t\rightarrow +\infty}\frac{1}{t}\log \# P(t)=
\lim_{t\rightarrow +\infty}\frac{1}{t}\log \# P_{reg}(t)=h_{top}(G),
\end{equation}
where $h_{top}(G)$ is the topological entropy of the geodesic flow $G$.

For compact rank one surfaces, K. Gelfert and B. Schapira (\cite{GS} Theorem 1.2) proved, among other results, a generalization of (1) to 
the case where geodesics are counted with weights. 
Namely, for H\"older continuous ``hyperbolic potentials'' $f$ 
(see \cite{GS} Theorem 1.2 and Lemma 3) which are constant 
on the singular set $sing$, they proved that there exists $\beta=\beta(f)> 0$ such that,
\begin{equation}
\limsup_{t\rightarrow +\infty}\frac{1}{t}\log
\sum_{\gamma \in P(t): \chi(\gamma)>\beta}
e^{\int_{0}^{l(\gamma)}f(g^{t}(\dot{\gamma} (0)))dt}=P(G;f),
\end{equation}
where $\chi(\gamma)$ is the smallest Lyapunov exponent of the periodic geodesic $\gamma$ and 
$P(G; f)$ the topological pressure of $G$ corresponding to the potential $f$. Recall briefly that the $P(G;f)$ satisfies the variational principle,
\[
P(G;f)=\sup \left (h_{m}(G)+\int fdm \right ),
\]
where the $\sup$ is over the set of invariant probability measures for $G$ and $h_{m}(G)=h_{m}(g^{1})$ is the measure theoretical entropy of the measure $m$ (see \cite{Hass}, \cite{Walters}).

Formula (2) suggests that the topological pressure of hyperbolic potentials can be recovered 
from the exponential growth rate of periodic geodesics 
with positive Lyapunov exponents (geodesics which stay far from the zero curvature part of the surface).
Note that the geodesics in the left hand side of (2) dont not intersect the singular set $sing$ 
(since in $sing$ the Lyapunov exponents are $0$). 
Thus the sum in (2)
is over the subset of regular geodesics where the curvature is negative. 
So, what is the contribution of the rest of all regular geodesics?
On the other hand, examples are given in \cite{GS} of non hyperbolic potentials which contradict (2).
Thus, what can be said if there is no hyperbolicity assumption on the potential $f$?
Note that  we always have,
\begin{equation}
\limsup_{t\rightarrow +\infty}\frac{1}{t}\log
\sum_{\gamma \in P(t)}
e^{\int_{0}^{l(\gamma)}f(g^{t}(\dot{\gamma} (0)))dt}\leq P(G; f),
\end{equation}
for any continuous function $f$ on $T^{1}M$.
Indeed, two geodesics representing two distinct free homotopy classes are $(\epsilon,t)$-separated
for any $\epsilon <inj(M)$ and $t>0$, i.e for any pair of geodesics $\gamma_{1}, \gamma_{2} \in P(t)$,
there exists $s\in [0,t]$ such that $d(\gamma_{1}(s), \gamma_{2}(s))>\epsilon$ \cite{kniper}.

Motivated by these questions, we prove (Theorem 1) that for H\"older 
continuous potentials (not necessarily hyperbolic) 
which are constant on the singular set $sing$ we have,
\begin{equation}
\limsup_{t\rightarrow +\infty}\frac{1}{t}\log
\sum_{\gamma \in P(t)}e^{\int_{\gamma}f}=
\limsup_{t\rightarrow +\infty}\frac{1}{t}\log
\sum_{\gamma  \in P_{reg}(t)}e^{\int_{\gamma}f}
=P(G; f),
\end{equation}
where we have set $\int_{\gamma}f:=\int
_{0}^{l(\gamma)}f(g^{t}(\dot{\gamma} (0)))dt$.
This means that, in spite of the existence of periodic regular geodesics 
which come close to the zero curvature part of the surface, the exponential growth 
is still determined by $P(f)$ as in the formula (2) of K. Gelfert and B. Schapira.
Thus, by $(2)$, $(3)$ and $(4)$, if no hyperbolicity is required on $f$, 
we only have an inequality in  (2) unless we add all regular geodesics with Lyapunov exponents close to zero. 

The proof of our result does not focus on large deviations technics, as in \cite{GS}. 
Indeed, large deviations estimates suppose that we are able to built a ``large
deviation functional'' on the space of invariant measures of the flow,
which must be positive on the complement of the subspace of equilibrium states. 
But this is not the case here due to the lack of hyperbolicity of the potential $f$.
Instead, we use a result of R. Bowen and D. Ruelle \cite{br}
to prove a formula for the topological pressure $P(G_{|\Lambda},f)$ of the restriction of the geodesic flow to a basic set
$\Lambda \subset T^{1}M$. To recover formula (4) 
we apply a result of K. Burns and K. Gelfert \cite{bg} which asserts that there exits
an increasing sequence $\Lambda_{1}\subset\Lambda_{2} \subset \cdots
\subset reg$
of basic sets, all contained in the regular part reg, such that 
the topological pressures $P(G|_{\Lambda_{l}},f)$ converge to $P(G; f)$ as $l\rightarrow +\infty$. In fact, they proved that
\[
\lim_{l\rightarrow +\infty}P(G|_{\Lambda_{l}},f)=
\sup_{l\geq 1}P(G|_{\Lambda_{l}},f)
=P(G; f).
\]
The main result of the paper is the following.
\begin{theorem} 
Let M be a compact connected rank 1 surface. 
Then, for any H\"older continuous potential $f$ on $T^{1}M$ 
which is constant on $sing$, the exponential growth rate of
rank 1 geodesics is given by the topological pressure, 
\begin{equation}
\limsup_{t\rightarrow +\infty}\frac{1}{t}
\log \sum_{\gamma \in P_{reg}(t)}e^{\int_{\gamma} f}=P(G; f).
\end{equation}
\end{theorem}
It is not clear how this result can be generalized, or not, to
compact rank 1 manifolds of $dim >2$.
Indeed, in the case of surfaces, 
closed hyperbolic subsets of $T^{1}M$  are contained 
in the regular set \cite{bg}, and this fact is 
not true in general (see \cite{cro, fich}). 

Each periodic geodesic $\gamma$ defines a probability 
measure on $T^{1}M$ with support in
$\{\dot{\gamma}(t), 0\leq t\leq l(\gamma)\}$ by,
\[
\delta_{\gamma}(A):=
leb (t\in [0, l(\gamma)]: g^{t}(\dot{\gamma} (0)) \in A)/
l(\gamma).
\]
The following corollary is an easy and standard consequence of Theorem 1 (on can adapt for example the method in \cite{am}). It gives a way of approximating equilibrium states.
\begin{corollary}
In the conditions of Theorem 1, the weak* limits 
of the measures
\[
\mu_{t}=\frac{1}{\sum_{\gamma \in P_{reg}(t)}e^{\int_{\gamma} f}}
\sum_{\gamma\in P_{reg}(t)}e^{\int_{\gamma} f}\delta_{\gamma}
\]
are measures of maximal pressure for $G$, i.e equilibrium
states of $G$ corresponding to the potential $f$.
\end{corollary}

As stated above, the main idea in our approach consists in
proving the main theorem above for the restriction of the geodesic flow to basic sets.
For this, recall that a basic set $\Lambda$ is a compact locally maximal
hyperbolic subset of $T^{1}M$ on which the flow is transitive. 
The hyperbolic set $\Lambda$ is closed and invariant with no fixed point.
The set $Per(\Lambda)$ of the periodic orbits of $G|_{\Lambda}$
is dense in $\Lambda$ and there is an open set $U \supset \Lambda$ with 
$\Lambda=\bigcap_{t\in \real} g^{t}U$.
By (\cite{bg} Lemma 4.1) the compact and invariant set $\Lambda$
is contained in the regular part $reg$ of $T^{1}M$.
In particular, the orbits in $Per(\Lambda)$ are regular.
Each $\gamma \in Per(\Lambda)$ is rank 1 and  
$\gamma$ is uniquely determined in its free homotopy class.
It is well known (see \cite{br} Theorem 3.3, \cite{hk} Corollary 20.3.8) that the geodesic 
flow on $\Lambda$ admits a unique equilibrium state $\mu_{f,\Lambda}$. 
The probability measure $\mu_{f,\Lambda}$ is ergodic and satisfies other interesting
properties that can be found in  \cite{br}. 
The following proposition gives an equidistribution result for the restriction 
of geodesic flow to the basic set $\Lambda$. Let $Per_{\Lambda}(t)$ be the subset of
$Per(\Lambda)$ of (primitive) periodic geodesics $\gamma$ with length $l(\gamma)\leq t$.
\begin{proposition}
Let $M$ be a compact finite dimensional smooth manifold 
of nonpositive curvature and $\Lambda \subset T^{1}M$ a basic set.
Let $f$ be a H\"older continuous function on $T^{1}M$.
\begin{enumerate}
\item The topological pressure of $G|_{\Lambda}$
is given by,
\begin{equation*}
P(G|_{\Lambda};f)=\lim_{t\rightarrow +\infty}\frac{1}{t}\log
\sum_{\gamma \in Per_{\Lambda}(t)}e^{\int_{\gamma}f}.
\end{equation*}
\item For any continuous function $K$ we have,
\begin{equation*}
\int K d\mu_{\Lambda,t}:=
\frac{\sum_{\gamma \in Per_{\Lambda}(t)}
e^{\int_{\gamma }f}\delta_{\gamma }(K)}
{\sum_{\gamma \in Per_{\Lambda}(t)}e^{\int_{\gamma}f}} \
\stackrel{t\rightarrow +\infty}{\longrightarrow} 
\int Kd\mu_{f,\Lambda}.
\end{equation*}
\end{enumerate}
\end{proposition}

\section{Proof of the main result}
\begin{proof}
Since the restriction of the potential $f$ to the singular part $sing$ is a constant function, by 
(\cite{bg} Theorem 1.4 and Theorem 6.1), 
we can find a family of compact locally maximal hyperbolic sets 
(basic sets) $\Lambda_{1} \subset \Lambda_{2} \cdots \subset reg$ such that,
\begin{equation}
\lim_{l\rightarrow +\infty}P(G|_{\Lambda_{l}};f)=\sup_{l\geq 1}P(G|_{\Lambda_{l}};f)
=P(G; f).
\end{equation}
For all $l\geq 1$ we have,
\begin{eqnarray*}
P(G|_{\Lambda_{l}};f) &=&
\lim_{t\rightarrow +\infty}\frac{1}{t}\log
\sum_{\gamma \in Per_{\Lambda_{l}}(t)}e^{\int_{\gamma}f} \ (by \ Proposition 1)\\
&\leq& \limsup_{t\rightarrow +\infty}\frac{1}{t}\log
\sum_{\gamma \in P_{reg}(t)}e^{\int_{\gamma}f} \ (since\ \Lambda_{l} \subset reg).
\end{eqnarray*}
By (6)
\begin{equation*}
P(G; f)\leq \limsup_{t\rightarrow +\infty}\frac{1}{t}\log
\sum_{\gamma \in P_{reg}(t)}e^{\int_{\gamma}f}.
\end{equation*}
Equality follows now from (3). This proves 
formula (5) and Theorem 1.
\end{proof}
\section{Proof of Proposition 1}
We begin with the following simple lemma.
\begin{lemma} Let $h$ be a smooth diffeomorphism acting
on a compact metric space $X$ with a unique equilibrium state $\nu_{f}$ corresponding to a continuous potential $f$.
Let $A\subset X$, $A \ne X$, be closed and invariant with $\nu_{f}(A)=0$. Then
\[
P(h|_{A};f) < P(h;f),
\]
where $P(h|_{A};f)$ and $P(h;f)$ are the topological pressure of the $h$ on $A$ and $X$ respectively.
\end{lemma}
\begin{remark} 
It is not clear if the conclusion of the above lemma remains true
if $\nu_{f}(A)\ne 0$.
However, if the potential $f$ is constant everywhere, then
under certain conditions on $h$, by (Lemma 6.2 \cite{bowen2}) we have $h_{top}(h|_{A})<h_{top}(h)$. 
\end{remark}
\begin{proof}
Let $\nu_{A}$ be an invariant probability measure for $h|_{A}$,
which satisfies $h(\nu_{A})+\int_{A} f d\nu_{A}=P(h|_{A};f)$.
$h$ being a diffeomorphism and $\nu_{A}(A)=1$, the measure $\nu_{A}$ is also an invariant measure for the dynamical system $(X,h)$. 
We have $\nu_{f}(A)=0$ so that $\nu_{A} \neq \nu_{f}$. 
Thus the lemma follows since $\nu_{f}$ is unique.
\end{proof}

A special flow (or suspension) is associated to the restriction
of the geodesic flow to the basic set $\Lambda$.
It was proved in \cite{bowen2} and (\cite{br}  Proposition 3.1)
that the geodesic flow $G$ on the basic set $\Lambda$
can be conjugated to a special flow $(\Lambda(\Sigma), \psi)$ 
over a topologically mixing subshift of finite type 
$\sigma :\Sigma \rightarrow \Sigma$ 
(the mixing property was proved in \cite{br}  Proposition 3.1). 
Furthermore, the conjugacy $\rho : \Lambda(\Sigma) \rightarrow
\Lambda$ is an isomorphism of the measurable
flows $(\Lambda(\Sigma), \psi)$ and  $(\Lambda,G)$.
In particular, the suspension admits a unique equilibrium
state $\mu_{f*,\Lambda}$ for the H\"older continuous potential 
on $\Lambda(\Sigma)$ given by $f^{*}=f\circ \rho$ and
$\mu_{f*,\Lambda}=\rho^{*}\mu_{f,\Lambda}$. 

Denote by $Per(\Lambda(\Sigma))$ the set of periodic orbits of the flow $\psi$ in $\Lambda(\Sigma)$ and, by $Per_{\Lambda(\Sigma)}(t)$ the subset of orbits $\tau$ with (least) period less than $t$. 
Recall that the suspension is expansive since this is the case for the subshift. The topological pressure of a suspension over a topologically mixing subshift of finite type is given by,
\begin{equation}
P(\psi; f^{*})=\lim_{t\rightarrow +\infty} 
\frac{1}{t}\log
\sum_{\tau \in Per_{\Lambda(\Sigma)}(t)}
e^{\int_{\tau}f^{*}}.
\end{equation}
We have $P(G|_{\Lambda}; f)=P(\psi; f^{*})$ (see \cite{br}).
There are two sets $Z\subset \Lambda$ and $Z' \subset \Lambda(\Sigma)$ 
such that $\rho: Z' \rightarrow Z$ is an homeomorphism \cite{bowen2}.
By (\cite{bowen2} Proposition 4.1),  
$\tau \in \Lambda(\Sigma)$ is periodic if and only if $\rho(\tau)$ is periodic. 
Thus,
\begin{equation}
\lim_{t\rightarrow +\infty}\frac{1}{t}\log
\sum_{\gamma \in Per_{\Lambda}(t): \dot{\gamma}(0)\in Z}e^{\int_{\gamma}f}=
\lim_{t\rightarrow +\infty}\frac{1}{t}\log
\sum_{\tau \in Per_{\Lambda(\Sigma)}(t):\tau(0) \in Z'}e^{\int_{\tau}f*}.
\end{equation}
There are two closed subsets (see \cite{bowen2} p452) $E^{s}, E^{u}$ 
which are invariant in one direction : 
$g^{t}(E^{s})\subset E^{s}$, $g^{-t}(E^{u})\subset E^{u}$
for $t\geq 0$, such that
\[
Per(\Lambda)=(Per(\Lambda)\cap Z) \cup (Per(\Lambda)\cap (E^{s} \cup E^{u})),
\]
and
\[
Per(\Lambda(\Sigma)=(Per(\Lambda(\Sigma))\cap Z') \cup 
(Per(\Lambda(\Sigma))\cap (\rho^{-1}(E^{s}) \cup \rho^{-1}(E^{u})).
\]
Set $E^{s,u}=E^{s} \cup E^{u}$. 
We have $\mu_{f^{*},\Lambda}(\rho^{-1}(E^{s,u}))=\mu_{f,\Lambda}(E^{s,u})=0$ (see \cite{bowen2}, \cite{br}).
Set $f^{*}_{1}(x):=\int_{0}^{1}f^{*}(\psi^{t}(x))dt$.
By Lemma 2 applied to 
$(\rho^{-1}(E^{s}), \psi^{1}, f^{*}_{1})$ and
$(\rho^{-1}(E^{u}), \psi^{-1}, f^{*}_{1})$, we have
\begin{equation}
\lim_{t\rightarrow +\infty}
\frac{\sum_{\tau \in Per_{\Lambda(\Sigma)}(t):\tau(0)\in
\rho^{-1}(E^{s,u})}
e^{\int_{\tau}f^{*}}}{\sum_{\tau \in Per_{\Lambda(\Sigma)}(t)}
e^{\int_{\tau}f^{*}}}=0.
\end{equation}
Indeed,
\[
P(\psi^{1}|_{\rho^{-1}(E^{s})};f^{*}_{1})< 
P(\psi^{1};f^{*}_{1})=P(\psi;f^{*}),
\]
\[
P(\psi^{-1}|_{\rho^{-1}(E^{u})};f^{*}_{1})< 
P(\psi^{-1};f^{*}_{1})=P(\psi;f^{*}).
\]
On the other hand, since the suspension is expansive, if two periodic points are not dynamically separated (i.e by $\psi$) then they belong
to the same periodic orbit of $\psi$. We have then,
\[
\limsup_{t\rightarrow +\infty}
\frac{1}{t} 
\log \sum_{\tau \in Per_{\Lambda(\Sigma)}(t): \tau(0)\in
\rho^{-1}(E^{s,u})}e^{\int_{\tau}f^{*}}\leq 
P(\psi^{1}|_{\rho^{-1}(E^{s})};f^{*}_{1})+P(\psi^{-1}|_{\rho^{-1}(E^{u})};f^{*}_{1}).
\]
Thus
\begin{equation}
\lim_{t\rightarrow +\infty}
\frac{\sum_{\tau \in Per_{\Lambda(\Sigma)}(t): \tau(0)\in Z'}
e^{\int_{\tau}f^{*}}}{\sum_{\tau \in Per_{\Lambda(\Sigma)}(t)}
e^{\int_{\tau}f^{*}}}=1,
\end{equation}
and then,
\begin{equation}
\lim_{t\rightarrow +\infty}\frac{1}{t}\log
\sum_{\tau \in Per_{\Lambda(\Sigma)}(t): \tau(0)\in Z'}e^{\int_{\tau}f^{*}}=
P(\psi;f^{*}).
\end{equation}
We deduce from (11), (8) and the fact that 
$P(\psi;f^{*} )=P(G|_{\Lambda}; f)$,
\begin{equation}
\lim_{t\rightarrow +\infty}\frac{1}{t}\log
\sum_{\gamma \in Per_{\Lambda}(t): \dot{\gamma}(0)\in Z}e^{\int_{\gamma}f}=P(G|_{\Lambda}; f).
\end{equation}
Finally, combining (12) and (3) we get, 
\[
\lim_{t\rightarrow + \infty}\frac{1}{t} \log 
\sum_{\gamma \in Per_{\Lambda}(t)}e^{\int_{\gamma}f}
=P(G|_{\Lambda}; f).
\]
This proves the first part of Proposition 1.

We prove now the second part of Proposition 1.
The following equidistribution result for the 
hyperbolic suspension
$(\Lambda(\Sigma), \psi)$ can be found in (\cite{polasterix} p120).
However, it can be easily proved using the asymptotic formula (7) and uniqueness of the equilibrium state $\mu_{f^{*}, \Lambda}$ for the suspension.
Recall that $\delta_{\tau}(K):=
\frac{1}{l(\tau)} \int_{0}^{l(\tau)} K(\tau (t))dt $.
\begin{theorem}[\cite{polasterix}]
For any continuous function $K:\Lambda(\Sigma) \rightarrow \real$ we have,
\begin{equation*}
\frac{\sum_{\tau \in Per_{\Lambda(\Sigma)}(t)}
e^{\int_{\tau}f^{*}}\delta_{\tau}(K)}
{\sum_{\tau \in Per_{\Lambda(\Sigma)}(t)}e^{\int_{\tau}f^{*}}} \
\stackrel{t\rightarrow +\infty}{\longrightarrow} \int Kd\mu_{f^{*},\Lambda}.
\end{equation*}
\end{theorem}
It follows from Theorem 2, (9) and (10),
\[
\mu_{\Lambda(\Sigma),t}^{Z'}(K):=
\frac{\sum_{\tau \in Per_{\Lambda(\Sigma)}(t):
\tau(0) \in Z'}
e^{\int_{\tau}f^{*}}\delta_{\tau}(K)}
{\sum_{\tau \in Per_{\Lambda(\Sigma)}(t):
\tau(0) \in Z'}e^{\int_{\tau}f^{*}}} \
\stackrel{t\rightarrow +\infty}{\longrightarrow} \int Kd\mu_{f^{*},\Lambda},
\]
for any continuous function $K$, i.e $\mu_{\Lambda(\Sigma),t}^{Z'}$
converges weakly to $\mu_{f^{*},\Lambda}$, as $t\rightarrow +\infty$.
Observe that we can not deduce this directly from Theorem 2 by
considering the function $1_{Z'}$; indeed, it is not clear
whether or not it's continuity (or equivalently the continuity of $1_{Z}$) can be deduced from \cite{bowen2} and the construction there.

Define the probability measures, 
\[
\mu_{\Lambda,t}^{Z}(K):=
\frac{\sum_{\gamma \in Per_{\Lambda}(t): \dot{\gamma}(0)\in Z}
e^{\int_{\gamma}f}\delta_{\gamma}(K)}
{\sum_{\gamma \in Per_{\Lambda}(t): \dot{\gamma}(0)\in Z}e^{\int_{\gamma}f}}.
\]
Since $\rho$ is an isomorphism, we have
$\rho^{*}\mu_{\Lambda,t}^{Z}=\mu_{\Lambda(\Sigma),t}^{Z'}$.
Now, $\rho$ being continuous, the measure $\mu_{\Lambda,t}^{Z}$ 
converges weakly to $\rho_{*}\mu_{f*,\Lambda}=\mu_{f,\Lambda}$. 

We show now how to deduce the convergence of the measures 
$\mu_{\Lambda,t}$ to the equilibrium state $\mu_{f,\Lambda}$, as
$t\rightarrow +\infty$. 
For any continuous function $K$ we have 
\begin{equation}
\mu_{\Lambda,t}(K)=\mu_{\Lambda,t}^{Z}(K) x_{t}+y_{t}(K)
\end{equation}
where,
\[
x_{t}=\frac{\sum_{\gamma \in Per_{\Lambda}(t): \dot{\gamma}(0)\in Z}
e^{\int_{\gamma}f}}{\sum_{\gamma \in Per_{\Lambda}(t)}
e^{\int_{\gamma}f}},
\]
and
\[
y_{t}(K)=\frac{\sum_{\gamma \in Per_{\Lambda}(t): \dot{\gamma}(0)\in Z^{c}}
e^{\int_{\gamma}f}\delta_{\gamma}(K)}{\sum_{\gamma \in Per_{\Lambda}(t)}
e^{\int_{\gamma}f}}.
\]
By part 1 of Proposition 1 and the same arguments in the above proof of this part, $x_{t}$ converges to 1 and $y_{t}$ converges to 0, as $t\rightarrow +\infty$. 
Finally, by (13) we see that 
$\mu_{\Lambda,t}(K)$ and $\mu_{\Lambda,t}^{Z}(K)$ converge to the same limit,
i.e $\mu_{f,\Lambda}(K)$. This proves that $\mu_{\Lambda,t}$
converges to $\mu_{f,\Lambda}$.

\end{document}